\documentclass[12pt]{article}
\usepackage{latexsym,amssymb,enumerate,float,geometry,cite}
\newtheorem{theorem}{Theorem}
\newtheorem{conjecture}[theorem]{Conjecture}
\newtheorem{claim}{Claim}

\begin{document}

%\linenumbers

\title{\Large A lower bound on the acyclic matching number of subcubic graphs}
\author{\large M. F\"{u}rst \and D. Rautenbach}
\date{}
\maketitle
\vspace{-10mm}
\begin{center}
{\small
Institute of Optimization and Operations Research, Ulm University,\\
Ulm, Germany, \texttt{maximilian.fuerst,dieter.rautenbach@uni-ulm.de}}
\end{center}

\begin{abstract}
The acyclic matching number of a graph $G$ is the largest size of an acyclic matching in $G$, that is, 
a matching $M$ in $G$ such that the subgraph of $G$ induced by the vertices incident to an edge in $M$ is a forest. 
We show that the acyclic matching number 
of a connected subcubic graph $G$ with $m$ edges 
is at least $m/6$ 
except for two small exceptions. 
\end{abstract}
{\small
\begin{tabular}{lp{12.5cm}}
\textbf{Keywords:} & 
Acyclic matching; subcubic graph
\end{tabular}
}

\section{Introduction}

We consider finite, simple, and undirected graphs, and use standard terminology and notation.
A matching $M$ in a graph $G$ is \textit{acyclic} \cite{gohehela} 
if the subgraph of $G$ induced by the set of vertices that are incident to some edge in $M$ is a forest,
and the \textit{acyclic matching number} $\nu_{ac}(G)$ of $G$
is the maximum size of an acyclic matching in $G$.
While the ordinary matching number $\nu(G)$ of $G$ is tractable \cite{ed},
it has been known for some time that the acyclic matching number is NP-hard for graphs of maximum degree $5$ 
\cite{gohehela,sp}. Recently, we \cite{fura} showed that just deciding the equality of $\nu(G)$ and $\nu_{ac}(G)$
is already NP-complete when restricted to bipartite graphs $G$ of maximum degree $4$.
The complexity of the acyclic matching number for cubic graphs is unknown.

In the present paper we establish a tight lower bound on the acyclic matching number of subcubic graphs.
Similar results were obtained for the matching number \cite{bideduflko,heye,owe,helora},
and also for the induced matching number \cite{hoqitr,jorasa,jo}.
Baste and Rautenbach \cite{bara} studied acyclic edge colorings,
and showed that the \textit{acyclic chromatic index} $\chi^\prime_{ac}(G)$ of a graph $G$, 
that is, the minimum number of acyclic matchings in $G$ into which the edge set of $G$ can be partitioned, is at most
$\Delta(G)^2$, where $\Delta(G)$ denotes the maximum degree of $G$. 
This implies $\nu_{ac}(G)\geq m(G)/\Delta(G)^2$, 
where $m(G)$ denotes the size of $G$,
which, for subcubic graphs, simplifies to $\nu_{ac}(G) \geq m(G)/9$.
This latter bound also follows from a lower bound \cite{jorasa} 
on the induced matching number,
which is always at most the acyclic matching number.
While the bound is tight for $K_{3,3}$,
excluding some small graphs allows a considerable improvement.
Let $K_4^+$ be the graph that arises by subdividing one edge of $K_4$ once. 

We prove the following.

\begin{theorem}\label{theorem1}
If $G$ is a connected subcubic graph that is not isomorphic to $K_4^+$ or $K_{3,3}$, 
then $\nu_{ac}(G) \geq m(G)/6$.
\end{theorem}
Since every subcubic graph $G$ of order $n(G)$ satisfies $m(G)\leq 3n(G)/2$, 
Theorem \ref{theorem1} is an immediate consequence of the following stronger result.
For two graphs $G$ and $H$, 
let $\kappa_G(H)$ denote the number of components of $G$ 
that are isomorphic to $H$. 
\begin{theorem}\label{theorem2}
If $G$ is a subcubic graph without isolated vertices, then 
\begin{eqnarray*}
\nu_{ac}(G) \geq \frac{1}{4}\left( n(G) - \kappa_G(K_{2,3}) - \kappa_G(K_4^+) - 2\kappa_G(K_{3,3}) \right).
\end{eqnarray*}
\end{theorem}
Note that Theorem \ref{theorem2} is tight;
examples are $K_4$, $K_{2,2}$, $K_{1,3}$, 
or the graph obtained from $K_{1,3}$ 
by replacing each endvertex with an endblock isomorphic to $K_{2,3}$. 
The proof of Theorem \ref{theorem2} is postponed 
to the second section.
The reduction arguments within that proof easily lead to a polynomial time algorithm 
computing acyclic matchings of the guaranteed size.

In a third section, we conclude with some open problems.

\section{Proof of Theorem \ref{theorem2}}\label{sectionproof2}

The proof is by contradiction.
Therefore, suppose that $G$ is a counterexample to Theorem \ref{theorem2} that is of minimum order $n$. 
A graph is {\it special} if it is isomorphic to $K_{2,3}$, $K_4^+$, or $K_{3,3}$.
Clearly, $G$ is connected, not special, and $n$ is at least $5$.
Note that $\nu_{ac}(G)<n/4$. 

We derive a contradiction using a series of claims.

\begin{claim} \label{c1}
No subgraph of $G$ is isomorphic to $K_4^+$.
\end{claim}
{\it Proof of Claim \ref{c1}:}
Suppose that $G$ has a subgraph $H$ that is isomorphic to $K_4^+$.
Let $v_1$, $v_2$, $v_3$, and $v_4$ be the vertices of degree $3$ in $H$, 
and let $u$ the vertex of degree $2$ in $H$.
Let $G'=G-\{v_1,v_2,v_3,v_4\}$. 
Since $G$ is connected, the graph $G'$ is connected.
Since $u$ has degree $1$ in $G'$, the graph $G'$ is not special.
By the choice of $G$, the graph $G'$ is no counterexample to Theorem \ref{theorem2},
and, hence, it has an acyclic matching $M'$ of size at least $n(G')/4=n/4-1$.
Adding the edge $v_1v_2$ to $M'$
yields an acyclic matching in $G$ of size at least $n/4$, 
which is a contradiction. 
$\hfill\Box$

\begin{claim} \label{c2}
No endblock of $G$ is isomorphic to $K_{2,3}$.
\end{claim}
{\it Proof of Claim \ref{c2}:}
Suppose that some endblock $B$ of $G$ is isomorphic to $K_{2,3}$.
Let $u$ be the unique cutvertex of $G$ in $B$.
Clearly, the vertex $u$ has degree $2$ in $B$. 
The graph $G' = G - (V(B)\setminus \{ u\})$ is connected,
and, since $u$ has degree $1$ in $G'$, it is not special.
Therefore, by the choice of $G$,
the graph $G'$ has an acyclic matching $M'$ of size at least $n(G')/4=n/4-1$.
Adding an edge of $B$ that is not incident to $u$ to $M'$
yields an acyclic matching in $G$ of size at least $n/4$, 
which is a contradiction. 
$\hfill\Box$

\begin{claim} \label{c3}
No two vertices of degree $1$ have a common neighbor.
\end{claim}
{\it Proof of Claim \ref{c3}:}
Suppose that $u$ and $v$ are two vertices of degree $1$, and that $w$ is their common neighbor.
Let $G' = G - \{u,v,w\}$.
Since $G'$ is connected and not isomorphic to $K_{3,3}$,
the choice of $G$ implies that $G'$ has an acyclic matching $M'$ of size at least $(n(G')-1)/4=n/4-1$.
Since $w$ does not lie on any cycle in $G$, adding the edge $uw$ to $M'$ 
yields an acyclic matching in $G$ of size at least $n/4$, 
which is a contradiction.
$\hfill\Box$

\begin{claim} \label{c4}
No vertex of degree $1$ is adjacent to a vertex that does not lie on a cycle.
\end{claim}
{\it Proof of Claim \ref{c4}:}
Suppose that $u$ is a vertex of degree $1$ that is adjacent to a vertex $v$ that does not lie on a cycle.
By Claim \ref{c3}, the graph $G' = G - \{u,v \}$ has no isolated vertex.
Since $G'$ has at most two components, and no component of $G'$ is isomorphic to $K_{3,3}$,
the choice of $G$ implies that $G'$ has an acyclic matching $M'$ of size at least $(n(G')-2)/4=n/4-1$.
Since $v$ does not lie on a cycle, adding the edge $uv$ to $M'$ 
yields an acyclic matching in $G$ of size at least $n/4$, 
which is a contradiction.
$\hfill\Box$
 
\begin{claim} \label{c5}
The minimum degree of $G$ is at least $2$.
\end{claim}
{\it Proof of Claim \ref{c5}:}
Suppose that $u$ is a vertex of degree $1$.
By Claim \ref{c4}, the neighbor $v$ of $u$ lies on a cycle $C$ in $G$.
Let $x$ and $w$ be the neighbors of $v$ on $C$.
 
First, suppose that $w$ has no neighbor of degree $1$.

If $G-\{ u,v,w\}$ contains an isolated vertex,
then this is necessarily the vertex $x$, and $N_G(x)=\{ v,w\}$.
In this case, let $G'=G-\{ u,v,w,x\}$.
Clearly, the graph $G'$ is connected and not isomorphic to $K_{3,3}$.
If isomorphic to $K_4^+$ or $K_{2,3}$, 
then it follows easily that $\nu_{ac}(G) \geq 3 > 9/4=n/4$,
which is a contradiction.
Hence, $G'$ is not special,
which implies that $G'$ has an acyclic matching $M'$ of size at least $n(G')/4=n/4-1$.
Adding the edge $uv$ to $M'$ yields an acyclic matching in $G$ of size at least $n/4$, 
which is a contradiction.
Hence, we may assume that $G'=G-\{ u,v,w\}$ has no isolated vertex.

Since there are at most three edges between $\{ u,v,w\}$ and $V(G')$ in $G$,
Claim \ref{c2} implies that at most one component of $G'$ is isomorphic to $K_{2,3}$.
By the choice of $G$,
this implies that $G'$ has an acyclic matching $M'$ of size at least $(n(G')-1)/4=n/4-1$.
Adding the edge $uv$ to $M'$ yields an acyclic matching in $G$ of size at least $n/4$, 
which is a contradiction.
Hence, by symmetry, we may assume that $x$ and $w$ both have a neighbor of degree $1$.

Let $y$ be a neighbor $w$ of degree $1$.
If $x$ and $w$ are adjacent, then $\nu_{ac}(G) = 2 >6/4=n/4$,
which is a contradiction.
Hence, $x$ and $w$ are not adjacent.
In view of the cycle $C$, the graph $G'=G - \{u,v,w,y\}$ is connected.
Since $G'$ has a vertex of degree $1$, it is not special,
which implies that $G'$ has an acyclic matching $M'$ of size at least $n(G')/4=n/4-1$.
Adding the edge $uv$ to $M'$ yields an acyclic matching in $G$ of size at least $n/4$, 
which is a contradiction.
$\hfill\Box$

\medskip

\noindent For a set $X$ of vertices of $G$, let $N_G[X]=\bigcup_{u\in X}N_G[u]$.

\begin{claim} \label{c6}
No subgraph of $G$ is isomorphic to $K_{2,3}$.
\end{claim}
{\it Proof of Claim \ref{c6}:}
Suppose that $G$ has a subgraph $H$ that is isomorphic to $K_{2,3}$. 
Claim \ref{c1} implies that $H$ is an induced subgraph of $G$.
Let $u_1$, $u_2$, and $u_3$ be the vertices of degree $2$ in $H$, 
and let $v_1$ and $v_2$ be the vertices of degree $3$ in $H$.
 
First, suppose that $u_1$ has degree $2$ in $G$.
Since $G$ is not special, we may assume that $u_2$ has degree $3$ in $G$.
By Claim \ref{c5}, the graph $G'=(V(H) \setminus \{ u_2\})$ has no isolated vertex,
and, since $u_2$ has degree $1$ in $G'$, it is not special. 
It follows that $G'$ has an acyclic matching $M'$ of size at least $n(G')/4=n/4-1$.
Adding the edge $u_1v_1$ to $M'$ yields an acyclic matching in $G$ of size at least $n/4$, 
which is a contradiction.
Hence, by symmetry, we may assume that all vertices in $U=\{ u_1,u_2,u_3\}$ have degree $3$ in $G$. 

Next, suppose that $u_1$ and $u_2$ have a common neighbor $u$ that is distinct from $v_1$ and $v_2$.
Let $G' = G - N_G[U]$. 
Note that there are at most $3$ edges between $N_G[U]$ and $V(G')$ in $G$.
By Claim \ref{c5}, the graph $G'$ has at most one isolated vertex, and,
by Claim \ref{c1}, at most one component of $G'$ is isomorphic to $K_{2,3}$.
Furthermore, the graph $G'$ does not have an isolated vertex as well as a component isomorphic to $K_{2,3}$.
This implies that $G'$ has an acyclic matching $M'$ of size at least $(n(G')-1)/4=n/4-2$.
Adding the two edges $uu_1$ and $u_3v_1$ to $M'$ yields an acyclic matching in $G$ of size at least $n/4$, 
which is a contradiction.
Hence, by symmetry, no two vertices in $U$ have a common neighbor that is distinct from $v_1$ and $v_2$.

The graph $G'$ that arises by contracting all edges of $H$
is simple and connected.
If $G'$ is special, then $G$ has order at most $11$, 
and an acyclic matching consisting of the three edges between $N_G[U]$ and $V(G)\setminus N_G[U]$ in $G$,
which is a contradiction.
Hence, $G'$ is not special, which implies that $G'$ has an acyclic matching $M'$ of size at least $n(G')/4=n/4-1$.
Let $M''$ be the acyclic matching in $G$ corresponding to $M'$.
Since $M''$ covers at most one vertex in $U$, say $u_1$, 
adding the edge $u_2v_1$ to $M''$ yields an acyclic matching in $G$ of size at least $n/4$, 
which is a contradiction.
$\hfill\Box$

\medskip

\noindent Claim \ref{c1}, Claim \ref{c6}, and the choice of $G$ imply that 
every proper induced subgraph $G'$ of $G$ with $i(G')$ isolated vertices 
has an acyclic matching $M'$ such that
\begin{eqnarray} \label{key}
|M'| \geq \frac{n(G')-i(G')}{4}.
\end{eqnarray}

\begin{claim} \label{c7}
No two vertices of degree $2$ are adjacent.
\end{claim}
{\it Proof of Claim \ref{c7}:}
Suppose that $u$ and $v$ are adjacent vertices of degree $2$, and that $w$ is the neighbor of $u$ distinct from $v$.
By Claim \ref{c5}, the graph $G' = G - \{u,v,w \}$ has at most one isolated vertex,
and, hence, by (\ref{key}), 
it has an acyclic matching $M'$ of size at least $(n(G')-1)/4=n/4-1$.
Adding the edge $uv$ to $M'$ yields a contradiction.
$\hfill\Box$

\begin{claim} \label{c8}
No vertex of degree $2$ lies on a triangle.
\end{claim}
{\it Proof of Claim \ref{c8}:}
Suppose that $u_1u_2u_3u_1$ is a triangle in $G$ such that $u_1$ has degree $2$.
By Claim \ref{c7}, the vertices $u_2$ and $u_3$ have degree $3$. 
Since $n\geq 5$, the graph $G' = G - \{u_1,u_2,u_3\}$ has no isolated vertex,
and, hence, by (\ref{key}), it has an acyclic matching $M'$ of size at least $n(G')/4>n/4-1$.
Adding the edge $u_1u_2$ to $M'$ yields a contradiction. 
$\hfill\Box$

\begin{claim} \label{c9}
No vertex of degree $2$ lies on a cycle of length $4$.
\end{claim}
{\it Proof of Claim \ref{c9}:}
Suppose that $u_1u_2u_3u_4u_1$ is a cycle in $G$ such that $u_1$ has degree $2$.
By Claims \ref{c7} and \ref{c8}, the vertices $u_2$ and $u_4$ have degree $3$, and are not adjacent.
By Claims \ref{c6} and \ref{c8}, the graph $G' = G - \{ u_1,u_2,u_3,u_4\}$ has no isolated vertex,
and, hence, by (\ref{key}), it has an acyclic matching $M'$ of size at least $n(G')/4=n/4-1$.
Adding the edge $u_1u_2$ to $M'$ yields a contradiction. 
$\hfill\Box$

\begin{claim} \label{c10}
No cycle of length $5$ contains two vertices of degree $2$.
\end{claim}
{\it Proof of Claim \ref{c10}:}
Suppose that the cycle $u_1u_2u_3u_4u_5u_1$ contains two vertices of degree $2$.
By Claim \ref{c7}, we may assume that $u_1$ and $u_4$ have degree $2$,
and that $u_2$, $u_3$, and $u_5$ have degree $3$.
Let $G'=G-(N_G[u_5]\cup \{ u_2,u_3\})$.
Since there are at most $4$ edges between $N_G[u_5]\cup \{ u_2,u_3\}$ and $V(G')$ in $G$,
the graph $G'$ has at most two isolated vertices,
and, hence, by (\ref{key}), it has an acyclic matching $M'$ of size at least $(n(G')-2)/4=n/4-2$.
Adding the edges $u_1u_2$ and $u_4u_5$ to $M'$ yields a contradiction. 
$\hfill\Box$

\begin{claim} \label{c11}
$G$ is cubic.
\end{claim}
{\it Proof of Claim \ref{c11}:}
Suppose that $u$ is a vertex of degree $2$.
By Claims \ref{c7}, \ref{c8}, and \ref{c9}, the neighbors of $u$, say $v$ and $w$, 
have degree $3$, are not adjacent, and have no common neighbor except for $u$.
Let $x$ be a neighbor of $v$ distinct from $u$.
By Claims \ref{c8}, \ref{c9}, and \ref{c10}, the graph $G' = G - \{u,v,w,x \}$ has no isolated vertex,
and, hence, by (\ref{key}), it has an acyclic matching $M'$ of size at least $n(G')/4=n/4-1$.
Adding the edge $uv$ to $M'$ yields a contradiction. 
$\hfill\Box$

\begin{claim} \label{c12}
$G$ is triangle-free.
\end{claim}
{\it Proof of Claim \ref{c12}:}
Suppose that $u_1u_2u_3u_1$ is a triangle in $G$.
By Claims \ref{c1} and \ref{c11}, the graph $G'=G-N_G[u_1]$ has no isolated vertex,
and, hence, by (\ref{key}), it has an acyclic matching $M'$ of size at least $n(G')/4=n/4-1$.
Adding the edge $u_1u_2$ to $M'$ yields a contradiction. 
$\hfill\Box$

\medskip

\noindent Let $C:u_1u_2\ldots u_gu_1$ be a shortest cycle in $G$.
For $i\in [g]$, let $v_i$ be the neighbor of $u_i$ not on $C$.
By Claim \ref{c12}, we have $g\geq 4$.

\begin{claim} \label{c13}
$g\geq 5$.
\end{claim}
{\it Proof of Claim \ref{c13}:}
Suppose that $g=4$.
By Claims \ref{c6} and \ref{c12}, the vertices $v_1$, $v_2$, $v_3$, and $v_4$ are distinct. 
Let $w_1$ and $w_2$ be the neighbors of $v_1$ distinct from $u_1$. 

First, suppose that $w_1=v_2$. 
By Claim \ref{c11}, the graph $G' = G - (N_G[v_1]\cup \{u_2,u_3,u_4 \})$ has at most one isolated vertex,
and, hence, by (\ref{key}), it has an acyclic matching $M'$ of size at least $(n(G')-1)/4=n/4-2$.
Adding the edges $u_1v_1$ and $u_2u_3$ to $M'$ yields a contradiction. 
Hence, we may assume, by symmetry, that $\{ v_1,v_2,v_3,v_4\}$ is independent.

Next, suppose that there is some vertex $x$ outside of $N_G[\{ v_1,u_1,u_3\}]$ 
such that $N_G(x)\subseteq  N_G[\{ v_1,u_1,u_3\}]$.
By Claim \ref{c6}, $x$ is not adjacent to both $u_2$ and $u_4$.
Hence, by Claim \ref{c11}, we may assume that $x$ is adjacent to $w_1$ but not to $u_2$.
By Claim \ref{c11}, the graph $G'=G-N_G[\{ v_1,u_1,u_3,w_1\}]$ has at most two isolated vertices,
and, hence, by (\ref{key}), it has an acyclic matching $M'$ of size at least $(n(G')-2)/4=n/4-3$.
Adding the edges $xw_1$, $u_1v_1$, and $u_2u_3$ to $M'$ yields a contradiction. 
Hence, we may assume that the graph $G'=G-N_G[\{ v_1,u_1,u_3\}]$ has no isolated vertex.
By (\ref{key}), $G'$ has an acyclic matching $M'$ of size at least $n(G')/4=n/4-2$.
Adding the edges $u_1v_1$ and $u_2u_3$ to $M'$ yields a contradiction. 
$\hfill\Box$

\begin{claim} \label{c14}
$g\geq 6$.
\end{claim}
{\it Proof of Claim \ref{c14}:}
Suppose that $g=5$.
By Claim \ref{c13}, the vertices $v_1$, $v_2$, $v_3$, $v_4$, and $v_5$ are distinct. 
Suppose that there is some vertex $x$ outside of $N_G[\{ u_1,u_2,u_4\}]$ 
such that $N_G(x)\subseteq  N_G[\{ u_1,u_2,u_4\}]$.
By Claims \ref{c11} and \ref{c13}, we obtain $N_G(x)=\{ v_1,v_2,v_4\}$.
By Claim \ref{c11}, the graph $G'=G-N_G[\{ v_1,u_1,u_2,u_4\}]$ has at most two isolated vertices,
and, hence, by (\ref{key}), it has an acyclic matching $M'$ of size at least $(n(G')-2)/4=n/4-3$.
Adding the edges $xv_1$, $u_1u_2$, and $u_3u_4$ to $M'$ yields a contradiction. 
Hence, we may assume that the graph $G'=G-N_G[\{ u_1,u_2,u_4\}]$ has no isolated vertex.
By (\ref{key}), the graph $G'$ has an acyclic matching $M'$ of size at least $n(G')/4=n/4-2$.
Adding the edges $u_1u_2$ and $u_3u_4$ to $M'$ yields a contradiction. 
$\hfill\Box$

\begin{claim} \label{c15}
$g\geq 7$.
\end{claim}
{\it Proof of Claim \ref{c15}:}
Suppose that $g=6$.
Let $w_1$ and $w_2$ be the neighbors of $v_1$ distinct from $u_1$.
By Claim \ref{c14}, the vertices $v_i$ for $i\in [6]\setminus \{ 4\}$, $w_1$, and $w_2$ are distinct. 
Suppose that there is some vertex $x$ outside of $N_G[\{ v_1,u_3,u_5,u_6\}]$ 
such that $N_G(x)\subseteq  N_G[\{ v_1,u_3,u_5,u_6\}]$.
By Claims \ref{c11} and \ref{c14}, we obtain that $x$ is adjacent to $v_3$,
to one vertex in $\{ v_5,v_6\}$, and 
to one vertex in $\{ w_1,w_2\}$.
Let $G'=G-N_G[\{ v_1,v_3,u_3,u_5,u_6\}]$.
By Claim \ref{c14}, no isolated vertex in $G'$ is adjacent to $u_2$ or $u_4$.
Since there are at most $10$ edges between $N_G[\{ v_1,v_3,u_3,u_5,u_6\}]$ and $V(G')$ in $G$,
this implies that $G'$ has at most two isolated vertices,
and, hence, by (\ref{key}), it has an acyclic matching $M'$ of size at least $(n(G')-2)/4=n/4-4$.
Adding the edges 
$xv_3$,
$u_1v_1$,
$u_2u_3$, and 
$u_5u_6$
to $M'$ yields a contradiction. 
Hence, we may assume that the graph $G'=G-N_G[\{ v_1,u_3,u_5,u_6\}]$ has no isolated vertex.
By (\ref{key}), the graph $G'$ has an acyclic matching $M'$ of size at least $n(G')/4=n/4-3$.
Adding the edges $u_1v_1$, $u_2u_3$, and $u_5u_6$ to $M'$ yields a contradiction. 
$\hfill\Box$

\medskip

\noindent We are now in a position to complete the proof.

First, suppose that $g$ is odd. 
If the graph $G'=G-N_G[\{ u_1,\ldots,u_{g-2}\}]$ has an isolated vertex,
then, by Claim \ref{c11}, there is a cycle of length at most $\left\lfloor\frac{g}{3}\right\rfloor+4$.
Since the last expression is less than $g$ for odd $g$ at least $7$,
it follows that $G'$ has no isolated vertex.
By (\ref{key}), the graph $G'$ has an acyclic matching $M'$ of size at least $n(G')/4=n/4-(g-1)/2$.
Adding the edges in $\{ u_{2i-1}u_{2i}:i\in [(g-1)/2]\}$ to $M'$ yields a contradiction. 
Hence, we may assume that $g$ is even.
Let $w_1$ and $w_2$ be the neighbors of $v_1$ distinct from $u_1$.
By the choice of $C$,
the vertices $v_i$ for $i\in [g]$, $w_1$, and $w_2$ are distinct. 
If the graph $G'=G-N_G[\{ v_1,u_1,\ldots,u_{g-2}\}]$ has an isolated vertex,
then, by Claim \ref{c11}, there is a cycle of length at most $\left\lfloor\frac{g}{3}\right\rfloor+5$.
Since the last expression is less than $g$ for even $g$ at least $8$,
it follows that $G'$ has no isolated vertex.
By (\ref{key}), the graph $G'$ has an acyclic matching $M'$ of size at least $n(G')/4=n/4-g/2$.
Adding the edges in $\{ u_1v_1\}\cup \{ u_{2i}u_{2i+1}:i\in [(g-2)/2]\}$ to $M'$ yields a contradiction,
which completes the proof.
$\hfill\Box$

\section{Conclusion}

We believe that Theorem \ref{theorem2} can be improved as follows.

\begin{conjecture} \label{conjecture1}
There is a constant $c$ such that 
$\nu_{ac}(G) \geq \frac{3n(G)}{11}-c$
for every connected subcubic graph $G$.
\end{conjecture}
Conjecture \ref{conjecture1} would be asymptotically best possible. 
If $H$ arises from a copy of $K_{1,2}$, 
where $u(H)$ denotes the vertex of degree $2$,
by replacing each endvertex with an endblock isomorphic to $K_{2,3}$,
and, for some positive integer $k$,
the connected subcubic graph $G_k$ 
arises from $k$ disjoint copies $H_1,\ldots,H_k$ of $H$
by adding, for every $i\in[k-1]$, 
an edge between $u(H_i)$ and some vertex of degree $2$ in $H_{i+1}$ 
that is distinct from $u(H_{i+1})$, 
then $\nu_{ac}(G_k)=3n(G_k)/11$.

For general maximum degree, we pose the following conjecture motivated by \cite{jo}.

\begin{conjecture}\label{conjecture2}
If $G$ is a graph of maximum degree $\Delta$ without isolated vertices,
then $$\nu_{ac}(G) \geq \min{ \left\{ \frac{2n(G)}{\left( \lceil \frac{\Delta}{2}\rceil + 1\right) \left( \lfloor \frac{\Delta}{2} \rfloor +1 \right)},
\frac{n(G)}{2\Delta} \right\} }.$$

\end{conjecture}
There should be better lower bounds on the acyclic matching number
for graphs of large girth, and methods from \cite{codara,flho,hera}
might be useful. Moreover, a lower bound as Conjecture \ref{conjecture2}, which is essentially tight for all possible
densities of a graph $G$ of bounded maximum degree, would be interesting, yet very challenging.

\end{document}